\documentclass[journal,twoside,web]{ieeecolor}
\usepackage{lcsys}
\usepackage{cite}
\usepackage{amsmath,amssymb,amsfonts}
\usepackage{graphicx}
\usepackage{textcomp}

\usepackage{mathtools}
\usepackage{multirow}
\usepackage{booktabs}
\usepackage{tabularx}
\usepackage{subcaption}
\usepackage{xcolor}  
\usepackage{svg}

\usepackage{algorithm}
\usepackage{algorithmicx}
\usepackage{algpseudocode}

\definecolor{purple}{rgb}{0.8,0,0.3}
\definecolor{deepblue}{rgb}{0.2,0,1}

\makeatletter
\let\NAT@parse\undefined
\makeatother
\usepackage{hyperref}
\hypersetup{
    colorlinks,
    linkcolor={deepblue},
    citecolor={purple},
    urlcolor={black}
}

\def\BibTeX{{\rm B\kern-.05em{\sc i\kern-.025em b}\kern-.08em
    T\kern-.1667em\lower.7ex\hbox{E}\kern-.125emX}}
\newtheorem{theorem}{Theorem}
\newtheorem{lemma}{Auxiliary Lemma}

\newtheorem{definition}{Definition}

\newtheorem{prooftheorem}{Proof of Theorem}
\newtheorem{prooflemma}{Proof of Auxiliary Lemma}

\begin{document}
\title{Optimal Control for Minimizing Inescapable Ellipsoids in Linear Periodically Time-Varying Systems under Bounded Disturbances}
\author{Egor Dogadin, Alexey Peregudin
\thanks{Egor Dogadin is with the Faculty of Control Systems and Robotics, ITMO University, St. Petersburg, Russia (email: egor.dogadin@icloud.com).}
\thanks{Alexey Peregudin is with the School of Electrical and Electronic Engineering, University of Sheffield, Sheffield, United Kingdom (e-mail: peregudin@sheffield.ac.uk).}
}

\maketitle
\thispagestyle{empty}

\begin{abstract}
This letter addresses optimal controller design for periodic linear time-varying systems under unknown-but-bounded disturbances. We introduce differential Lyapunov-type equations to describe time-varying inescapable ellipsoids and define an integral-based measure of their size. To minimize this measure, we develop a differential Riccati equation-based approach that provides exact solutions for state-feedback, observer synthesis, and output-feedback control. A key component is a systematic procedure for determining the optimal time-varying parameter, reducing an infinite-dimensional optimization to a simple iterative process. A numerical example validates the method's effectiveness.
\end{abstract}

\begin{IEEEkeywords}
Invariant ellipsoids, time-varying systems, optimal control, optimal filtering, output feedback.
%Enter key words or phrases in alphabetical order, separated by commas. For a list of suggested keywords, send a blank e-mail to keywords@ieee.org or visit \underline{http://www.ieee.org/organizations/pubs/ani\_prod/keywrd98.txt}
\end{IEEEkeywords}

\section{Introduction}
\label{sec:introduction}
\IEEEPARstart{O}{ptimal} control under unknown-but-bounded disturbances remains a challenge in control theory.  While LQG, $\mathcal{H}_2$, and $\mathcal{H}_\infty$ controllers perform well under specific disturbances (e.g., white noise or harmonic inputs), they are generally suboptimal when disturbances are only known to be bounded, as in practical systems such as autonomous drones, marine vehicles, and robotic manipulators. In such cases, external influences like wind gusts, water currents, or actuator nonlinearities often lack a precise structure but can be bounded in magnitude. This has motivated the use of \emph{invariant (inescapable) ellipsoids} to characterize reachability bounds.

The concept of inescapable ellipsoids under bounded inputs dates back to \cite{Schweppe1973, Boyd1994}, which established fundamental equations forming the basis for controller design. Early applications, such as \cite{Abedor1994}, used these ellipsoids for peak-to-peak minimization, but \cite{Venkatesh1995} showed that minimizing the ellipsoid alone does not guarantee the smallest peak-to-peak gain, leading to alternative performance measures.

A different approach emerged with the introduction of a new trace-based optimality criterion \cite{Khlebnikov-1}, making it conceptually similar to $\mathcal{H}_2$ control (see \cite{Doyle1989}), and opening possibilities for controller design based on Linear Matrix Inequalities (LMIs).  This framework was extended to state-feedback \cite{Khlebnikov-1}, output-feedback \cite{Khlebnikov-2}, and observer design \cite{Khlebnikov-3}, and applied to nonlinear extensions \cite{e1,x7}, robust control \cite{x6}, and sliding-mode control \cite{e3,y3}, with a comprehensive summary in \cite{Poznyak2014}.  These formulations relied heavily on LMIs, often introducing robustness along with some degree of conservatism.

Exact solutions to the problems in \cite{Khlebnikov-1, Khlebnikov-2, Khlebnikov-3} were recently derived using Riccati equations \cite{peregudin2023, peregudin2024}. By introducing dual ellipsoids and duality relations, these works provided the first exact solution to the inescapable ellipsoid minimization problem for continuous-time LTI systems with output-feedback control. Similar ideas were later extended to discrete-time systems \cite{dogadin2024optimal}. However, minimization of inescapable ellipsoids for LTV systems remains an open problem, as existing studies \cite{Garcia2019, Azhmyakov2019} remain conservative and lack exact solutions.

To bridge this gap, this letter extends the Riccati-based approach in \cite{peregudin2024} to periodic LTV systems. Specifically, we introduce an optimality measure aligned with LTV $\mathcal{H}_2$ control, using a time-averaged trace-based criterion (see \cite{LTV}).

A key challenge in this framework is selecting the optimal parameter $\alpha$ (or its time-varying counterpart $\alpha(t)$) for constructing the inescapable ellipsoid \cite{Abedor1994, Khlebnikov-1, peregudin2024}. Even in LTI systems, choosing $\alpha$ is non-trivial, and in LTV systems, where $\alpha$ varies over time, the problem becomes even more complex \cite{Schweppe1973}. Previous approaches often assumed a constant $\alpha$, limiting optimality \cite{Poznyak2014}. We propose a method to compute the optimal $\alpha(t)$, making the problem tractable for LTV systems and also offering a more direct approach for the LTI case.

\vspace{0.5\baselineskip}

\textbf{Contributions of this Letter:}
\begin{itemize} 
 \item Exact (necessary and sufficient) characterization of periodic time-varying inescapable ellipsoids.
    \item An approach for periodic LTV systems that provides \textit{optimal} state-feedback, observer, and output-feedback controller under bounded disturbances.
    \item A systematic method to find the optimal $\alpha(t)$, reducing an infinite-dimensional search to a simple iteration.
\end{itemize}

\section{Analysis: Time-Varying Inescapable Ellipsoids}
Consider a linear periodically time-varying system
\begin{equation} \label{System}
    \begin{cases}
        \dot{x}(t) = A(t)x(t) + B(t)w(t), \\
        z(t) = C(t)x(t),
    \end{cases}
\end{equation}
where $x(t) \in \mathbb{R}^n$, $w(t) \in \mathbb{R}^m$, $z(t) \in \mathbb{R}^k$, and $A(\cdot)$, $B(\cdot)$, $C(\cdot)$ are real, piecewise-continuous, $T$-periodic matrix-valued functions, for some $T > 0$ satisfying $A(t) = A(t+T)$, $B(t) = B(t+T)$, $C(t) = C(t+T)$. Assume $(A(\cdot),B(\cdot))$ is controllable and $(C(\cdot),A(\cdot))$ is observable.

Suppose the system is subject to unit-bounded inputs 
\begin{equation} \label{Input_bound}
    \forall t \in \mathbb{R} \; : \; \lVert w(t) \rVert \le 1,
\end{equation} 
where $\|\cdot\|$ denotes the Euclidean norm. This unit bound incurs no loss of generality. In the general LTV case, it is natural to consider disturbances confined within a time-varying ellipsoid \( w(t)^\top W(t) w(t) \le 1 \), with \( w_{\min} I \preceq W(t) \preceq w_{\max} I \), as in \cite{Schweppe1973}. In periodic settings, \( W(t) \) is typically periodic. Such bounds can be naturally absorbed into the input matrix by the rescaling  \( B_{\text{new}}(t) = B(t)(W(t))^{-1/2} \), recovering \eqref{Input_bound}.

\begin{definition}[(Periodic inescapable ellipsoids)]
    Periodic inescapable ellipsoids $\mathcal{E}_P(\cdot)$ of the system \eqref{System}--\eqref{Input_bound} are defined as the time-varying sets
\begin{equation*}
    \mathcal{E}_P(t) \vcentcolon = \left\{\, {\rm x}\in \mathbb{R}^n \; \mid \; {\rm x}^\top P(t)^{-1} { \rm x} \le 1 \, \right\},
\end{equation*}
where $P(t)=P(t+T) \succ 0$, and the following inescapability property is satisfied:
\begin{equation} \label{Definition_Inescapability}
    x(t_0) \in \mathcal{E}_P(t_0) \;\; \Rightarrow \;\; \forall t \ge t_0 \; : \; x(t) \in \mathcal{E}_P(t).
\end{equation} 
\iffalse
\begin{equation*}
    \mathcal{E}_{P'}(t) \text{ s.t. \eqref{Definition_Periodicity}, \eqref{Definition_Inescapability}}, \;
        \mathcal{E}_{P'}(t) \subseteq \mathcal{E}_P(t)
    \; \Rightarrow \; \mathcal{E}_{P'}(t) = \mathcal{E}_P(t).
\end{equation*}  
\fi \vspace{-\baselineskip}
\end{definition}

\iffalse
Time-varying ellipsoids have been considered in control problems before (see \cite{Schweppe1973, Balandin2020}), but mostly in different settings and under conditions different from \eqref{Input_bound}.
\fi

For LTI systems, ellipsoids inescapable under \eqref{Input_bound} are known as invariant (attractive) ellipsoids \cite{Khlebnikov-1,Khlebnikov-2,Khlebnikov-3,e1,x7,x6,e3,y3,Poznyak2014} and are described by a constant matrix $P$. In this case, they are time-independent, and the periodicity condition is obsolete. Moreover, a stationary ellipsoid $\mathcal{E}_P$ is inescapable for \eqref{System}--\eqref{Input_bound} with constant $A$, $B$, $C$  iff there exists a constant $\alpha>0$ such that
% \begin{equation*}
%      \mathcal{L}_{P}\, \Big\{ A + \frac{\alpha}{2}I, \, \frac{1}{\alpha} B B^\top \Big\} =0,
% \end{equation*}
\begin{equation*}
      AP + PA^\top + \alpha P + \frac{1}{\alpha} B B^\top \preceq \, 0, \quad P \succ 0.
\end{equation*}
For a fixed $\alpha$, the smallest (inclusion-minimal) ellipsoid is obtained by replacing $\preceq$ with $=$. Hence, in the LTI case, the family of minimal inescapable ellipsoids is parametrized by the single parameter $\alpha$.

The following theorem characterizes the periodic inescapable ellipsoids of an LTV system. 

\begin{theorem} \label{Analysis_equation}
$\mathcal{E}_P(\cdot)$ is a periodic inescapable ellipsoid for the system \eqref{System}--\eqref{Input_bound} if and only if the matrix-valued function $P(t)=P(t+T)\succ 0$ satisfies the differential inequality
% \begin{equation}\label{P_equation}
%      \dot{P}(t) = \mathcal{L}_{P(t)}\, \Big\{ A(t) + \frac{\alpha(t)}{2}I, \, \frac{1}{\alpha(t)} B(t) B(t)^\top \Big\}
% \end{equation}
\begin{equation*}
    \phantom{1} \hspace{-0.2cm} \dot{P}(t) \!  \,\succeq\, \! A(t)P(t) + P(t)A(t)^{\! \top} \! \! + \alpha(t) P(t)  +  \frac{1}{\alpha(t)} B(t) B(t)^{\! \top} \hspace{-0.36cm} \phantom{1}
\end{equation*}
 for some positive function $\alpha(t)>0$. For a fixed $\alpha(\cdot)$, the smallest (inclusion-minimal) ellipsoid is obtained from
 \begin{equation}\label{P_equation}
    \phantom{1} \hspace{-0.2cm} \dot{P}(t) \! = \! A(t)P(t) + P(t)A(t)^{\! \top} \! \! + \alpha(t) P(t)  +  \frac{1}{\alpha(t)} B(t) B(t)^{\! \top} \hspace{-0.36cm}\phantom{1}.
\end{equation}
\vspace{-0.5\baselineskip}
\end{theorem}

Each minimal periodic inescapable ellipsoid is uniquely determined by \(\alpha(\cdot)\) and may be denoted \(\mathcal{E}_{\alpha}(\cdot) \doteq \mathcal{E}_{P}(\cdot) \), where $P$ is given by \eqref{P_equation}. By periodic differential Lyapunov theory \cite{PeriodicLyapunov}, the ellipsoid is well-defined if \(\alpha(\cdot)\) is \(T\)-periodic, \(A(t) + \frac{\alpha(t)}{2}I\) is exponentially stable, and \((A(\cdot), B(\cdot))\) is controllable.

Equation \eqref{P_equation} dates back to \cite{Schweppe1973}, which also inspired the underlying idea of Theorem~\ref{Analysis_equation}, though it was not stated there in its present form. We provide a simpler, self-contained proof of Theorem \ref{Analysis_equation} (included in the \hyperlink{appendix}{Appendix} along with other proofs). 
%Notably, although this time-varying equation appeared in \cite{Schweppe1973}, 
Notably, later works \cite{Abedor1994, Khlebnikov-1, Khlebnikov-2, Khlebnikov-3, e1, x7, x6, e3, y3, Poznyak2014} used only a constant \(\alpha\) even for LTV systems \cite{Garcia2019, Azhmyakov2019}, thus restricting the class of candidate ellipsoids and potentially missing the optimal one. This letter reconsiders \(\alpha\) as time-dependent and introduces an effective method to compute the optimal \(\alpha(\cdot)\), minimizing the periodic inescapable ellipsoid. All subsequent results presented in the letter are, to the best of our knowledge, entirely new.

We introduce the following measure of the ellipsoid's size:
\begin{equation} \label{Size}
    \operatorname{size} (\mathcal{E}_{\alpha}) \doteq \frac{1}{T} \int_{0}^{T} \mathrm{trace} \left( C(t) P(t)C(t)^{\!\top} \right) dt,
\end{equation}
which can be interpreted as the average sum of squares of the semiaxes of the image of the ellipsoid $\mathcal{E}_{\alpha}(t)$ under the mapping $x \mapsto z=C(t)x$. This measure aligns with established theory; for instance, replacing \( P(t) \) in~\eqref{Size} with the controllability Gramian of system~\eqref{System} yields the squared \( \mathcal{H}_2 \)-norm of the system (see~\cite{LTV}). The following duality result holds:

\begin{theorem} \label{equivalence_of_norms}
Let $P(t)=P(t+T)\succ 0$ be given by \eqref{P_equation}, and let $Q(t)=Q(t+T)\succ 0$ satisfy the dual equation
% \begin{equation}\label{Q_equation}
%      -\dot{Q}(t) = \mathcal{L}_{Q(t)}\, \Big\{ A(t)^{\! \top} + \frac{\alpha(t)}{2}I, \, C(t)^{\! \top} C(t) \Big\}.
% \end{equation}
\begin{equation}\label{Q_equation}
     -\dot{Q}(t)  = Q(t)A(t) + A(t)^{\! \top} Q(t) + \alpha(t) Q(t)  +  C(t)^{\! \top} C(t).
\end{equation}
The size of the ellipsoid $\mathcal{E}_{\alpha}(\cdot)$ can then be expressed in terms of the dual matrix $Q(t)$ as follows:
\begin{equation*}
    \operatorname{size} (\mathcal{E}_{\alpha}) = \frac{1}{T} \int_{0}^{T} \frac{\mathrm{trace} \left( B(t)^{\! \top} Q(t) B(t)\right)}{\alpha(t)} dt.
\end{equation*}
\vspace{-\baselineskip}
\end{theorem}

The system's sensitivity to disturbance \( w \) in \eqref{System} is given by the \textit{smallest} inescapable ellipsoid over all \( \alpha(\cdot) \). The following theorem provides a method to determine this tightest bound.

\begin{theorem} \label{Optimal_alpha_theorem}
    (I) The optimal time-dependent parameter $\alpha^*(t) = \arg \min_{\alpha} \operatorname{size} (\mathcal{E}_{\alpha}) $ satisfies
    \begin{equation} \label{Optimal_alpha}
        \alpha^*(t)   =  \sqrt{\frac{\mathrm{trace}(B(t)^{\!\top}Q(t) B(t))}{\mathrm{trace}(Q(t)P(t))}},
    \end{equation}
where $P(t)$, $Q(t)$ are such that \eqref{P_equation}, \eqref{Q_equation} with $\alpha(t)=\alpha^*(t)$. \linebreak
(II) With any $\alpha_0(t)=\alpha_0(t+T)>0$, the fixed-point iteration 
\begin{equation} \label{Fixed_point}
        \alpha_{i+1}(t)   =  \sqrt{\frac{\mathrm{trace}(B(t)^{\!\top}Q_i(t) B(t))}{\mathrm{trace}(Q_i(t)P_i(t))}}, 
    \end{equation}
where $P_i(t)$, $Q_i(t)$ are such that \eqref{P_equation}, \eqref{Q_equation} with $\alpha(t)=\alpha_i(t)$, converges to the globally optimal $\alpha^*(t)$.
\end{theorem}

\section{Design: Minimizing Inescapable Ellipsoids}

\subsection{Optimal State-Feedback Controller}

Consider a linear periodically time-varying system
\begin{equation} \label{State-Feedback_System}
    \begin{cases}
        \dot{x}(t) = A(t)x(t) + B(t)u(t) + B_w(t)w(t), \\
        z(t) = C(t)x(t)+D(t)u(t),
    \end{cases}
\end{equation}
where $x(t) $ is the state, $z(t) $ is the regulated output, $u(t) $ is the control, $w(t) $ is the disturbance and $A $, $B $, $C $, $D$, $B_w$ are piecewise-continuous $T$-periodic matrix-functions. Here  $(A(\cdot),B(\cdot))$ is stabilizable, $(C(\cdot),A(\cdot))$ is observable, and $D(t)^{\!\top}\!\begin{bmatrix}
        C(t) & D(t)
    \end{bmatrix}= \begin{bmatrix}
        0 & I
    \end{bmatrix}$.
The latter standard assumption (see \cite{Doyle1989}) simplifies formulas but can be omitted at the cost of more complex controller equations.

Consider a linear time-varying controller:
\begin{equation} \label{State-Feedback_Controller}
u(t) = K(t)x(t),
\end{equation}
the corresponding closed-loop system \eqref{Input_bound}, \eqref{State-Feedback_System}-\eqref{State-Feedback_Controller}, and the family of periodic inescapable ellipsoids $\mathcal{E}_{\alpha, K}(\cdot)$ associated with it. In line with \eqref{Size}, the size of each ellipsoid is given by 
\[
\operatorname{size} (\mathcal{E}_{\alpha, K}) = \frac{1}{T} \int_{0}^{T} \mathrm{trace} \left( C_K(t) P(t)C_K^\top(t) \right) dt,
\]
where $C_K(t)=C(t)+D(t)K(t)$, and $P$ is obtained from the closed-loop version of \eqref{P_equation} and depends on both $\alpha$ and $K$.

The state-feedback design objectives are as follows: \hypertarget{S1}{} \hypertarget{S2}{}
\begin{enumerate}
    \item[(S1)] \textit{Fixed-$\alpha$ optimality}. For a given $\alpha(\cdot)$, find the optimal controller gain $ K^*_\alpha(t) = \arg \min_{K} \operatorname{size} (\mathcal{E}_{\alpha, K})$. 
    \item[(S2)] \textit{Global optimality}. Find the optimal $\alpha(\cdot)$ and controller gain pair $ (\alpha^*(t), K^*(t)) = \arg \min_{\alpha, K} \operatorname{size} (\mathcal{E}_{\alpha, K})$. 
\end{enumerate}

The following theorem provides a method to obtain optimal controllers that fulfill the design objective (\hyperlink{S1}{S1}).

\begin{theorem} \label{Theorem_KQ} \color{black} 
For a given $\alpha(\cdot)>0$, the optimal controller gain $ K(t) = K^*_\alpha(t) = \arg \min_{K} \operatorname{size} (\mathcal{E}_{\alpha, K})$ is given by
\begin{equation} \label{State-Feedback_K}
    K(t) = - B (t)^{\!\top} Q(t),
\end{equation}
where $Q(t)=Q(t+T) \succ 0$ is obtained from
% \begin{equation} \label{State-Feedback_Q}
%     \dot{Q}(t)  \! =\! - \mathcal{R}_{ Q(t)} \Big\{A(t)^{\!\top}  \!  +  \frac{\alpha(t)}{2},  C(t)^{\!\top} \! C(t), B(t)B(t)^{\!\top} \! \Big\}.
% \end{equation}
\begin{equation} \label{State-Feedback_Q}
\begin{aligned}
    -\dot{Q}(t) & =  Q(t)A(t) + A(t)^\top Q(t) + \alpha(t) Q(t) \\
     & + C(t)^\top C(t) - Q(t)B(t)B(t)^\top Q(t).
    \end{aligned}
\end{equation}
% { \color{orange}
% The associated ellipsoid size attains its minimum:
%     \begin{equation*}
%         \min \operatorname{size} (\mathcal{E}_{\alpha, K}) =\frac{1}{T} \int_{0}^{T} \frac{1}{\alpha(t)} \mathrm{trace} (B_w^\top(t)Q(t)B_w(t)) dt.
%     \end{equation*}
%     }
% (II) The globally optimal controller gain \(K(t) = K^*(t)\), which minimizes $\operatorname{size} (\mathcal{E}_{\alpha, K})$ over all $\alpha(\cdot)$, is given by \eqref{State-Feedback_K}, where $Q(t)$, $P(t)$, $\alpha(t) = \alpha^*(t)$ satisfy \eqref{State-Feedback_Q} and 
% \begin{equation*}
% \begin{aligned}
% & \alpha^2(t) \mathrm{trace}(Q(t)P(t))  =  \mathrm{trace}(B_w^\top(t)Q(t) B_w(t)), \\
% &{\dot{P}(t)  =  \mathcal{L}_{P(t)} \! \Big\{A_K(t) \! +\! \frac{\alpha(t)}{2},  \frac{B_w(t)B_w(t)^\top}{\alpha(t)}\Big\}}, 
% \end{aligned}
% \end{equation*}
% where $A_K(t)= A(t)+B(t)K(t)$ and $P(t)=P(t+T) \succ 0$.
 \vspace{-0.5\baselineskip}
\end{theorem}
\color{black}

Theorem \ref{Theorem_KQ} provides the exact solution to problem (\hyperlink{S1}{S1}). The solution to problem (\hyperlink{S2}{S2}) is discussed in Section \ref{Iterative_section}.

\vspace{5px}

% \color{orange}

% \begin{theorem}
%     The locally optimal controller gain $K(t)$ that minimizes the $\varepsilon$-norm can be iteratively obtained by the Algorithm \ref{State-Feedback_Algorithm}.
% \end{theorem}

% \begin{algorithm}
%     \caption{Gain $K(t)$ $\varepsilon$-optimization} \label{State-Feedback_Algorithm}
%     \begin{algorithmic}[1]
%     \Require Initial guess: $\alpha_0 (t) = \alpha_0 (t+T) \in \mathbb{R}_+$, $\epsilon\in \mathbb{R}_+$.
%     \For{$i = 0, 1, 2, \cdots$}
%        \State Find $Q_i(t)= \Ric_{\alpha_i(t)}^* \big( A(t), B(t), C(t) \big)$.
%        \vspace*{4px}
%     \State Compute $K_i(t) =  - B^\top (t) Q_i(t)$. 
%     \vspace*{4px}
%     \State Find $P_i(t)=\Lyap_{\alpha_i(t)} \big(A(t) + B(t)K_i(t),B_w(t)\big)$.
%     \vspace*{4px}
%     \State \begin{minipage}[t]{0.93\linewidth} Use matrices $Q_i(t)$ and $P_i(t)$ to determine the next $\alpha_{i+1}(t)$: \end{minipage}
%     \[
%         \alpha_{i+1}(t) = \sqrt{\frac{\mathrm{trace}(B_w^\top(t)Q_i(t) B_w(t))}{\mathrm{trace}(Q_i(t)P_i(t))}}.
%     \]
%        \vspace*{1px}
%     \If{$\| \alpha_{i+1}-\alpha_i\|<\epsilon$}
%     \vspace*{4px}
%         \State Break out of the loop.
%     \vspace*{4px}
%     \EndIf
%     \EndFor
%     %\EndProcedure
% \end{algorithmic}
% \end{algorithm}

% \color{black}

\subsection{Optimal Observer (Filter)}

Consider a linear periodically time-varying system:
\begin{equation} \label{Filtering_System}
%\mathcal{S}_W \, : \
\begin{cases}
    \dot x(t) = A(t) x(t) + B(t) w(t), \\
    y(t) = C(t) x(t) + D(t)w(t),
\end{cases}
\end{equation}
where $x(t)$ is the state, $y(t)$ is the measured output, $w(t)$ is the disturbance, $A$, $B$, $C$, $D $ are  piecewise-continuous $T$-periodic matrix-functions,  $(A(\cdot),B(\cdot))$ is controllable, $(C(\cdot),A(\cdot))$ is detectable, and $D(t) \begin{bmatrix}
        B (t)^{\!\top}  & D (t)^{\!\top}
    \end{bmatrix}  = \begin{bmatrix}
       0  &  I
    \end{bmatrix}  $. The latter standard assumption (see \cite{Doyle1989}) simplifies formulas but can be omitted at the cost of more complex observer equations.

The system is connected to an observer:
\begin{equation} \label{Filtering_Observer}
    %\mathcal{S}_L \, : \
    \begin{cases}
        \dot {\hat x}(t) = A(t) \hat x(t) + L(t) (C(t)\hat x(t) - y(t)) , \\
        z(t) = C_w(t) (x(t) - \hat x(t)),
    \end{cases}
\end{equation}
where $\hat x(t) $ is the estimated state, $z(t) $ is the regulated output (``weighted''  error $e(t) = x(t) - \hat x (t)$), and $C_w$ is piecewise-continuous and $T$-periodic. The error dynamics are given by
\begin{equation} \label{Filtering_error}
   \begin{cases}
    \dot e(t) = A_L(t)  e(t) + B_L(t) w(t), \\
    z(t) =  C_w(t) e(t) ,
\end{cases}
\end{equation}
where $A_L(t) = A(t) + L(t)C(t)$, $B_L(t) = B(t) + L(t) D(t)$.

Note that the system \eqref{Filtering_error} is in the form \eqref{System} and consider the family of periodic inescapable ellipsoids $\mathcal{E}_{\alpha, L}(\cdot)$ associated with it. For each ellipsoid, the size $\operatorname{size} (\mathcal{E}_{\alpha, L})$ is given by \eqref{Size}, where $C=C_w$ and $P$ depends on both $\alpha$ and $L$.

The observer (filter) design objectives are as follows: \hypertarget{F1}{} \hypertarget{F2}{}
\begin{enumerate}
    \item[(F1)] \textit{Fixed-$\alpha$ optimality}. For a given $\alpha(\cdot)$, find the optimal observer gain $ L^*_\alpha(t) = \arg \min_{L} \operatorname{size} (\mathcal{E}_{\alpha, L})$. 
    \item[(F2)] \textit{Global optimality}. Find the optimal $\alpha(\cdot)$ and observer gain pair $ (\alpha^*(t), L^*(t)) = \arg \min_{\alpha, L} \operatorname{size} (\mathcal{E}_{\alpha, L})$. 
\end{enumerate}

The following theorem provides a method to obtain optimal observers that fulfill the design objective (\hyperlink{F1}{F1}).

\begin{theorem} \label{Theorem_LP} \color{black}
For a given $\alpha(\cdot)>0$, the optimal observer gain $ L(t) = L^*_\alpha(t) = \arg \min_{L} \operatorname{size} (\mathcal{E}_{\alpha, L})$ is given by
\begin{equation} \label{Filtering_L}
    L(t) = - \alpha(t) P(t) C (t)^{\!\top},
\end{equation}
where $P(t)=P(t+T) \succ 0$ is obtained from
% \begin{equation} \label{Filtering_P}
%     \dot{P}(t)  \! =\!  \mathcal{R}_{ P(t)}  \! \Big\{A(t) \! + \! \frac{\alpha(t)}{2}, \frac{B(t)B(t)^{\! \top}\!}{\alpha(t)} \!,  \alpha(t)C(t)^{\! \top}\!  C(t) \Big\}.
% \end{equation}
\begin{equation} \label{Filtering_P}
\begin{aligned}
    \dot{P}(t) &=  A(t)P(t) + P(t) A(t)^\top  + \alpha(t) P(t) \\
     &+ \frac{1}{\alpha(t)}B(t) B(t)^\top - \alpha(t) P(t)C(t)^\top C(t) P(t).
    \end{aligned}
\end{equation}
% { \color{orange}
% The associated ellipsoid size attains its minimum:
%     \begin{equation*}
%         \min \operatorname{size} (\mathcal{E}_{\alpha, L}) =\frac{1}{T} \int_{0}^{T} \frac{1}{\alpha(t)} \mathrm{trace} (C_w(t)P(t)C_w^\top(t)) dt.
%     \end{equation*}
%     }
% (II) The globally optimal observer gain \(L(t) = L^*(t)\), which minimizes $\operatorname{size} (\mathcal{E}_{\alpha, L})$ over all $\alpha(\cdot)$, is given by \eqref{Filtering_L}, where $Q(t)=Q(t+T) \succ 0$, $P(t)$, $\alpha(t) = \alpha^*(t)$ satisfy \eqref{Filtering_P} and 
% \begin{equation*}
% \begin{aligned}
% & \alpha^2(t) \mathrm{trace}(Q(t)P(t))  =  \mathrm{trace}(B_L^\top(t)Q(t)  B_L(t)), \\
% & \dot{Q}(t)   = - \mathcal{L}_{ Q(t)} \Big\{A_L(t)^{\!\top}    +  \frac{\alpha(t)}{2},  C_w(t)^{\!\top} \! C_w(t)\! \Big\}. 
% \end{aligned}
% \end{equation*}
 \vspace{-0.5\baselineskip}
\end{theorem}
\color{black}

Theorem \ref{Theorem_LP} provides the exact solution to problem (\hyperlink{F1}{F1}). The solution to problem (\hyperlink{F2}{F2}) is discussed in Section \ref{Iterative_section}.

\vspace{5px}

% \color{orange}

% \begin{theorem}
%     The locally optimal observer gain $L(t)$ that minimizes the $\varepsilon$-norm can be iteratively obtained by the Algorithm \ref{Filtering_Algorithm}.
% \end{theorem}

% \begin{algorithm}
%     \caption{Gain $L(t)$ $\varepsilon$-optimization} \label{Filtering_Algorithm}
%     \begin{algorithmic}[1]
%     \Require Initial guess: $\alpha_0 (t) = \alpha_0 (t+T) \in \mathbb{R}_+$, $\epsilon \in \mathbb{R}_+ $.
%     \For{$i = 0, 1, 2, \cdots$}
%        \State Find $P_i(t) = \Ric_{\alpha_i(t)} \big( A(t), B(t), C(t) \big)$. 
%        \vspace*{4px}
%     \State Compute $L_i(t) = - \alpha_i(t) P_i(t) C^\top (t)$.
%        \vspace*{4px}
%     \State Find $Q_i(t)= \Lyap_{\alpha_i (t)}^* \big(A(t) + L_i(t)C(t), C_w(t) \big) $.
%        \vspace*{4px}
%     \State \begin{minipage}[t]{0.93\linewidth} Use matrices $P_i(t)$ and $Q_i(t)$ to compute the next $\alpha_{i+1}(t)$:
%     \[
%         \alpha_{i+1}(t) = \sqrt{\frac{\mathrm{trace}(\tilde B_i^\top(t) Q_i(t) \tilde B_i(t))}{\mathrm{trace}(Q_i(t)P_i(t))}},
%     \] 
%     \vspace*{4px}
%     where $\tilde B_i(t) = \left(B(t) +L_i(t) D(t) \right)$.\end{minipage}
%        \vspace*{4px}
%     \If{$\| \alpha_{i+1}-\alpha_i\|<\epsilon$}
%     \vspace*{4px}
%         \State Break out of the loop.
%     \vspace*{4px}
%     \EndIf
%     \EndFor
%     %\EndProcedure
%     \end{algorithmic}
%     \end{algorithm}

% \color{black}

\subsection{Optimal Output-Feedback Controller} \label{Section III-C}

Consider a periodically time-varying system:
\begin{equation} \label{Output-Feedback_system}
\begin{cases}
    \dot x(t) = A(t) x(t) + B_2(t) u(t) + B_1(t) w(t), \\
    y(t) = C_1(t) x(t) + D_1(t)w(t),\\
    z(t) = C_2(t) x(t) + D_2 (t) u(t),
\end{cases}
\end{equation}
where $x(t)$ is the state, $y(t) $ is the measured output, $z(t) $ is the regulated output, $u(t)$ is the control, $w(t)$ is the disturbance, $\hat x(t) $ is the estimated state, $A$, $B_i$, $C_i$, $D_i$, are piecewise-continuous $T$-periodic matrix-functions, $(A(\cdot),B_1(\cdot))$ is controllable, $(C_2(\cdot),A(\cdot))$ is observable, $(A(\cdot),B_2(\cdot))$ is stabilizable, $(C_1(\cdot),A(\cdot))$ is detectable,  $D_1(t) \begin{bmatrix}
    B_1(t)^{\!\top} \! & D_1(t)^{\!\top}  
\end{bmatrix} \! = \! \begin{bmatrix}
    0 & I
\end{bmatrix}$,  $D_2(t)^{\!\top} \! \begin{bmatrix}
    C_2(t) & D_2(t)
\end{bmatrix}  \!  = \!\begin{bmatrix}
    0 & I
\end{bmatrix}$.

Consider a linear time-varying controller:
\begin{equation} \label{Output-Feedback_controller}
    \begin{cases} \!
        \dot {\hat x}(t) \! = \! A(t) \hat x(t)\! + \!B_2(t)u(t)\! + \! L(t) (C_1(t)\hat x(t) \!- \! y(t)), \hspace{-0.65cm} \phantom{1}  \\
       \! u(t) \! = \! K(t)\hat x,
    \end{cases}
\end{equation}
the corresponding closed-loop system \eqref{Input_bound}, \eqref{Output-Feedback_system}-\eqref{Output-Feedback_controller}, and the family of periodic inescapable ellipsoids $\mathcal{E}_{\alpha, K,L}(\cdot)$ associated with it. The size $\operatorname{size} (\mathcal{E}_{\alpha, K, L})$ of each ellipsoid is given by \eqref{Size}, where $C$ and $P$ correspond to the closed-loop system and depend on the choice of $\alpha$, $K$, and $L$.

The output-feedback design objectives are as follows:  \hypertarget{O1}{} \hypertarget{O2}{}
\begin{enumerate}
    \item[(O1)] \textit{Fixed-$\alpha$ optimality}. For a given $\alpha(\cdot)$, find the controller gains $ (K^*_\alpha(t),L^*_\alpha(t)) = \arg \min_{K,L} \operatorname{size} (\mathcal{E}_{\alpha, K,L})$. 
    \item[(O2)] \textit{Global optimality}. Find the optimal $\alpha$ and controller gains $ (\alpha^*\!(t), K^*\!(t), L^*\!(t)) \! = \! \arg \! \min_{\alpha, K,L} \operatorname{size} (\mathcal{E}_{\alpha,  K,L})$. 
\end{enumerate}

The following theorem provides a method to obtain optimal controllers that fulfill the design objective (\hyperlink{O1}{O1}).

\begin{theorem} \label{Theorem_KQLP} \color{black} 
For a given $\alpha(\cdot)>0$, the optimal controller $ (K(t), L(t)) = (K^*_\alpha(t),L^*_\alpha(t)) = \arg \min_{K,L} \operatorname{size} (\mathcal{E}_{\alpha, K,L})$ is
\begin{equation} \label{Output-Feedback}
\left\{ \begin{aligned}
    &\eqref{State-Feedback_K} \, \\
    &\eqref{State-Feedback_Q} \,
\end{aligned}\right| \begin{aligned}
    & B = B_2, \\
    & C = C_2,
\end{aligned} \quad \quad
\left\{ \begin{aligned}
    &\eqref{Filtering_L} \, \\
    &\eqref{Filtering_P} \,
\end{aligned}\right| \begin{aligned}
    & B = B_1, \\
    & C = C_1.
\end{aligned}
\end{equation}
% { \color{orange}
% The associated minimum ellipsoid size $ \min \operatorname{size} (\mathcal{E}_{\alpha, K, L})$ is:
%   \begin{equation*}
%     \begin{aligned}
%         &\int_{0}^{T}  \frac{\mathrm{trace} \big(B_1^\top \! (t)Q(t)B_1(t) + \alpha(t)K(t) P(t) K^\top \! (t)\big)}{\alpha (t) T} dt \\
%         =&\int_{0}^{T}  \frac{\mathrm{trace} \big(\alpha(t)C_2 (t)P(t)C_2^\top \!(t) + L^\top \! (t)Q(t) L(t)\big)}{\alpha (t) T} dt.
%         \end{aligned}
%     \end{equation*}
%     \vspace{0\baselineskip}
%     }
 \vspace{-0.5\baselineskip}
\end{theorem}

Theorem \ref{Theorem_KQLP} gives the exact solution to problem (\hyperlink{O1}{O1}), and for this problem, the separation principle holds: optimal gains $K_\alpha^*$ and $L_\alpha^*$ can be obtained independently for a fixed $\alpha(\cdot)$.  

The solution to problem (\hyperlink{O2}{O2}) is discussed in Section \ref{Iterative_section}. It is worth noting that, for this problem, the separation principle no longer holds because the optimal $\alpha^*$ that solves (\hyperlink{O2}{O2}) is generally different from ones that solve (\hyperlink{S2}{S2}) and (\hyperlink{F2}{F2}).

\vspace{5px}

\subsection{Iterative Approach Towards Global Optimality} \label{Iterative_section}

Minimizing inescapable ellipsoids is efficiently achieved via LMIs \cite{Khlebnikov-3, Poznyak2014} or Riccati equations \cite{peregudin2024, dogadin2024optimal} when $\alpha$ is fixed, but optimizing over all $\alpha$ remains challenging. For LTI systems, one can evaluate a range of numerical $\alpha$ values to find the optimum. In the LTV case, $\alpha(\cdot)$ is a function, making the search space infinite-dimensional.

Using Theorem \ref{Optimal_alpha_theorem}, we propose an efficient algorithm that seeks $\alpha^*(\cdot)$, addressing problems (\hyperlink{S2}{S2}), (\hyperlink{F2}{F2}) and (\hyperlink{O2}{O2}).

\iffalse
\textbf{Algorithm for Controller Optimization Over $\alpha$.} \hypertarget{Algorithm}{}
\begin{enumerate}
    \item Choose the initial guess $\alpha_0(t)=\alpha_0(t+T)>0$.
    \item Solve (\hyperlink{S1}{S1}), (\hyperlink{F1}{F1}) or (\hyperlink{O1}{O1}) to find $K(t)$, $L(t)$, or both.
    \item Using step 2 results, obtain the system in the form \eqref{System} with the closed-loop matrices $A(t)$, $B(t)$ and $C(t)$.
    \item With \eqref{P_equation}, \eqref{Q_equation}, find $P(t)$, $Q(t)$ of the closed-loop system.
    \item Apply the update rule \eqref{Fixed_point} to find $\alpha_{i+1}(t)$.
    \item Repeat steps 2–5 until convergence to obtain a \textcolor{blue}{candidate}~$\alpha^*(t)$.
    \item Compute the corresponding $K^*(t)$, $L^*(t)$, or both. 
\end{enumerate}
\fi

\begin{algorithm}[H] 
\caption{Controller Optimization Over $\alpha$}
\label{algorithm}
\begin{algorithmic}[1]
\State Choose the initial guess $\alpha_0(t) = \alpha_0(t+T) > 0$.
\State Solve (\hyperlink{S1}{S1}), (\hyperlink{F1}{F1}) or (\hyperlink{O1}{O1}) to find $K(t)$, $L(t)$, or both.
\State Obtain the system in the form \eqref{System} with closed-loop matrices $A(t)$, $B(t)$ and $C(t)$.
\State Use \eqref{P_equation}, \eqref{Q_equation} to find $P(t)$, $Q(t)$ of the closed-loop system.
\State Apply the update rule \eqref{Fixed_point} to compute $\alpha_{i+1}(t)$.
\State Repeat steps 2–5 until convergence to a candidate $\alpha^*(t)$.
\State Compute the corresponding $K^*(t)$, $L^*(t)$, or both.
\end{algorithmic}
\end{algorithm}

Theorem \ref{Optimal_alpha_theorem} guarantees that the update rule \eqref{Fixed_point} leads to a minimum; for the analysis problem, this minimum is global owing to the convexity of $\alpha \mapsto \operatorname{size}(\mathcal{E}_\alpha)$. However, for the closed-loop design, it is known \cite[Remark 3]{peregudin2023} that the function $\alpha \mapsto \operatorname{size}(\mathcal{E}_{\alpha, K, L})$ may be non-convex. Consequently, there is no immediate proof of the algorithm’s convergence in the general case. In practice, it is typically sufficient to employ constant initial choices $\alpha_0(t)\equiv c$, with $c>0$, to obtain the global optimum with Algorithm \ref{algorithm}.

\section{Numerical Example} \label{Section:Example}

\begin{figure}[t]
    \centering
    %\vspace{0.08cm}
    \includegraphics[width = 0.41\textwidth]{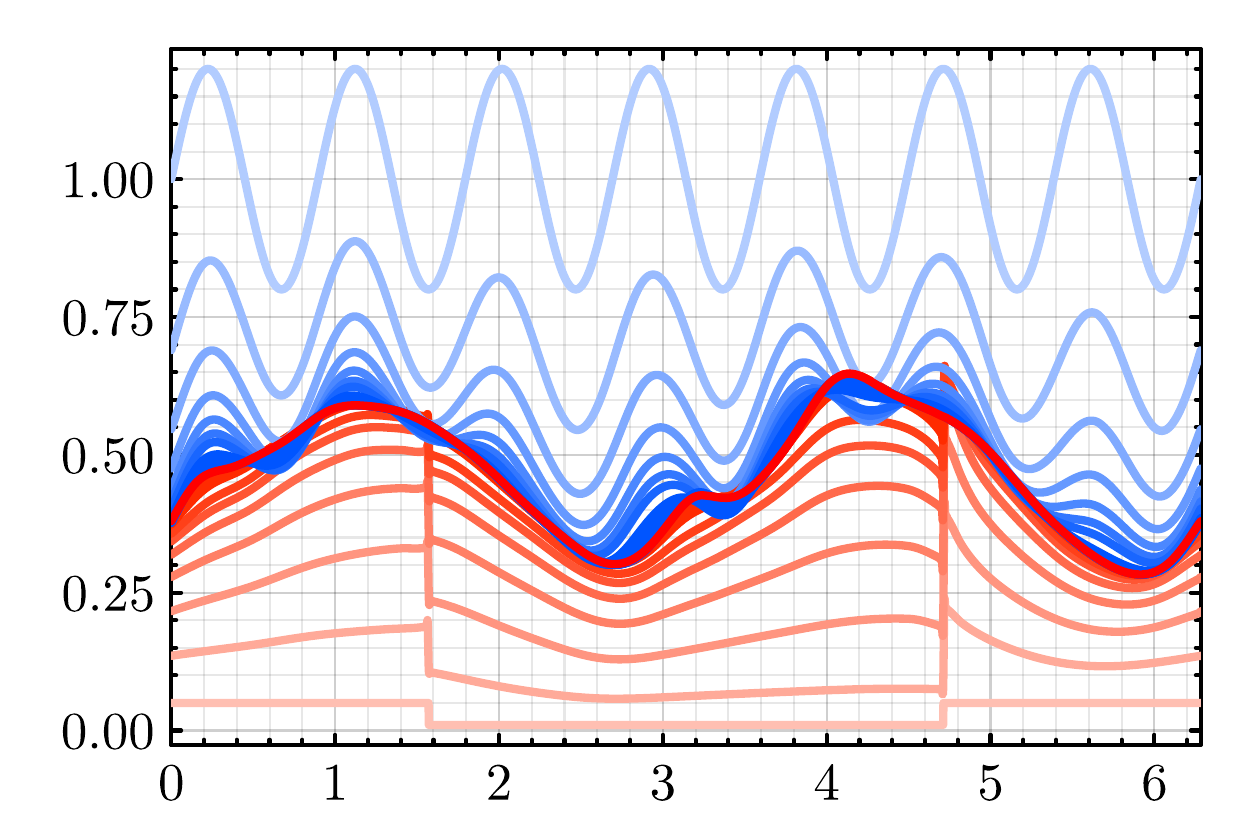} \hspace*{-0.4cm}
   \raisebox{0.6cm}{\includegraphics[width = 0.06\textwidth]{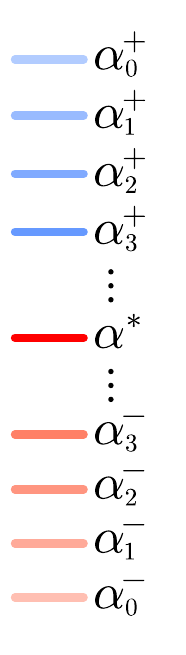}
   \hspace{-0.5cm}}
   \makebox[0pt][l]{\hspace{-8.7cm}
        \rotatebox{90}{\begin{minipage}{4.3cm}
            \[ \qquad  \alpha(t) \]
        \end{minipage}} \
    }
   \\
    \vspace*{-0.8cm}\[t,\mathrm{s} \ \ \] \\
    \vspace*{-0.3cm}
    %[width = 0.48\textwidth]
    \caption{Time-varying parameters $\alpha_i^+(t)$ and $\alpha_i^-(t)$ for the example in Section \ref{Section:Example} at different iterations, both converging to the optimal function $\alpha^*(t)$ from two different initial guesses.}
    \label{fig:alpha}
\end{figure}

Consider a system \eqref{Output-Feedback_system} with matrices
\begin{equation*}
    A \! = \! \! \begin{bmatrix}
        0 & 0 & 1 & 0 \\
        0 & 0 & 0 & 1 \\
        -\frac{k(t)}{J_1} & \frac{k(t)}{J_1} & -\frac{\mu(t)}{J_1} & \frac{\mu(t)}{J_1}\\
        \frac{k(t)}{J_2} & -\frac{k(t)}{J_2} & \frac{\mu(t)}{J_2} & -\frac{\mu(t)}{J_2}
    \end{bmatrix}\!\!, \,   B_1 \!  = \! \!\begin{bmatrix}
    0 & 0 & 0\\
    0 & 0 & 0 \\
    \frac{b_1(t)}{J_1} & 0 & 0 \\
    \frac{b_2(t)}{J_2} & 0 & 0
    \end{bmatrix}\!\!,
\end{equation*}
\begin{equation*}
     B_2 \! = \! \! \begin{bmatrix}
        0 \\ 0 \\ \frac{1}{J_1} \\ 0
    \end{bmatrix}\! , \; C_1 = \begin{bmatrix}
        1 & 0 \\
        0 & 1 \\
        0 & 0 \\
        0 & 0
    \end{bmatrix}^{\! \top}\!\!\!\!, 
    \; C_2 \! = \!\!\begin{bmatrix}
        c_1(t) & 0 \\
        c_2(t) & 0 \\
        0 & 0 \\
        0 & 0
    \end{bmatrix}^{\! \top} \!\!\!\!, \;
    D_1 \! = \!\! \begin{bmatrix}
        0 & 0 \\
        1 & 0 \\
        0 & 1 
    \end{bmatrix}^{\! \top}\!\!\!\!,
\end{equation*}
and $D_2 = \begin{bmatrix}
    0 & 1
\end{bmatrix}^{\! \top}$. This system models the dynamics of a coupled two-body system with time-varying stiffness $k(t)$ and damping $\mu(t)$. A common example is a gear pair with periodic mesh stiffness \cite{machines}, though its applications extend far beyond.

For this example, the following parameters were chosen: $k(t) = k_0(1 + 0.2\sin(t) +0.3\cos(2t))$,
 $\mu(t) = \mu_0(1+ 0.2\sin(t) - 0.3\cos(2t))$, $b_1 (t) = \frac{1}{3}\sin^2(t)$, $b_2(t) =  \frac{1}{10}\cos^2(t) $, $ c_1(t) =  20\sin^2(t)$, $c_2(t) = 7.5\cos^2(t)$, $J_1=0.5 \ \mathrm{kg\cdot m^2}$, $J_2=0.8 J_1$, $k_0=10 \ \mathrm{N\cdot m}$, $\mu_0=\frac{1}{300} \ \mathrm{N\cdot m\cdot s}$, $T=2\pi \ \mathrm{s}$. 

\begin{figure}[t!]
    \centering
    %\vspace{0.08cm}
    \hspace{-8cm}
    \includegraphics[width = 0.41\textwidth]{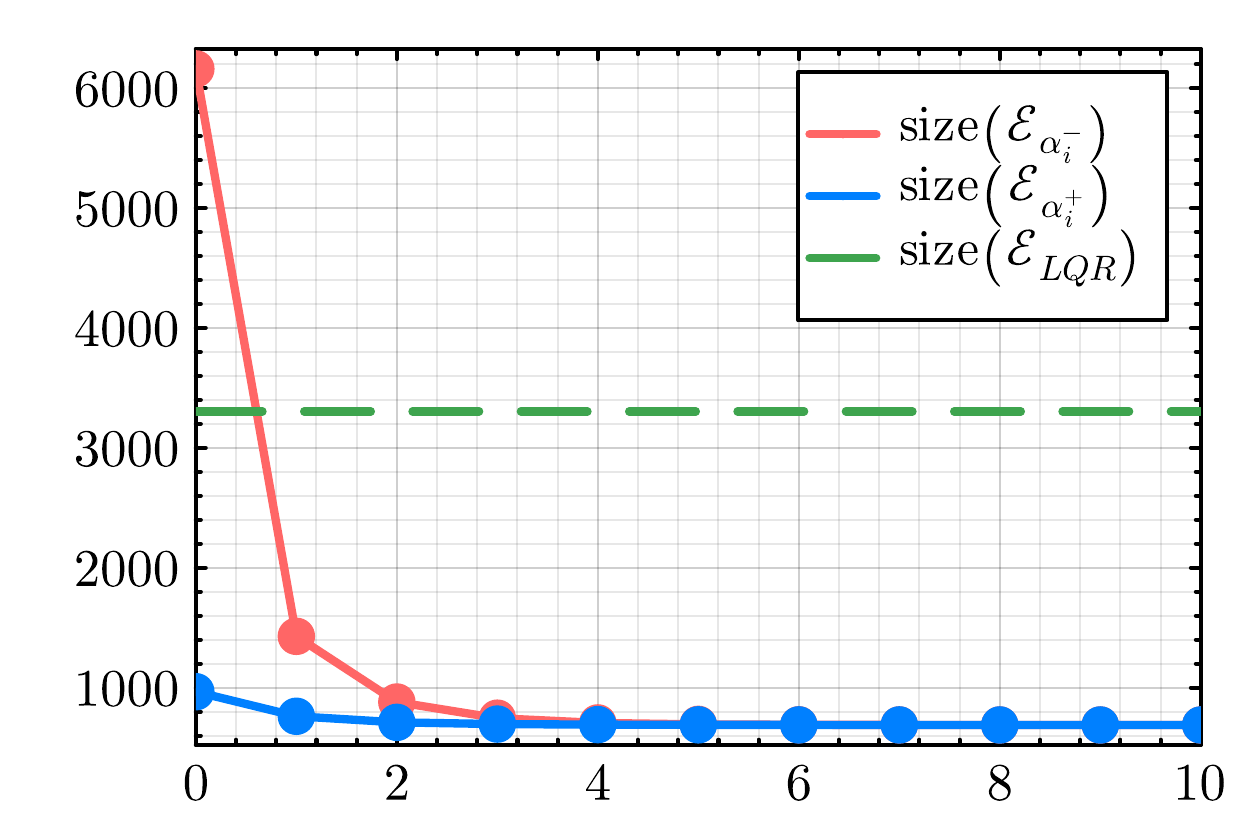}
    \hspace{-8cm}
    \rotatebox{90}{\begin{minipage}{5cm}
        \[ \ \mathrm{size}(\, \mathcal{E} \,) \]
    \end{minipage}} \\
    \vspace*{-0.2cm}\hspace{0.2cm}Iteration number
    %[width = 0.48\textwidth]
    \caption{Inescapable ellipsoid sizes in Section \ref{Section:Example} as a function of iteration number $i$, along with the inescapable ellipsoid size for the baseline LQR-Kalman filter combination.}
    \label{fig:size}
\end{figure}

The optimal controller was obtained using \eqref{Output-Feedback} and the algorithm in Section \ref{Iterative_section}. Figure \ref{fig:alpha} illustrates the convergence of $\alpha_i(t)$ to $\alpha^*(t)$ from two initial guesses: 
 $\alpha_0^-(t) = 0.01 $ if $t \mod T \in[\frac{\pi}{2}, \frac{3\pi}{2}]$ or $0.05$ otherwise and $\alpha_0^+(t) = 1 + \sin(7t)/5$, both reaching the same limit. Figure \ref{fig:size} shows the norm minimization during this process.

For comparison, we selected the LQR and Kalman filter combination. In line with the system description~\eqref{Output-Feedback_system}, the LQR weights were set as \( Q = C_2^\top C_2 \), \( R = 1 \); for the Kalman filter, \( Q = B_1 B_1^\top \), \( R = I \). Applying this controller, we computed the minimal inescapable ellipsoid size, shown in Figure \ref{fig:size}, which is larger than that of the proposed controller. 

\begin{figure*}[t!]
\centering
\begin{subfigure}{0.32\textwidth}
\centering
\includegraphics[width=0.98\textwidth]{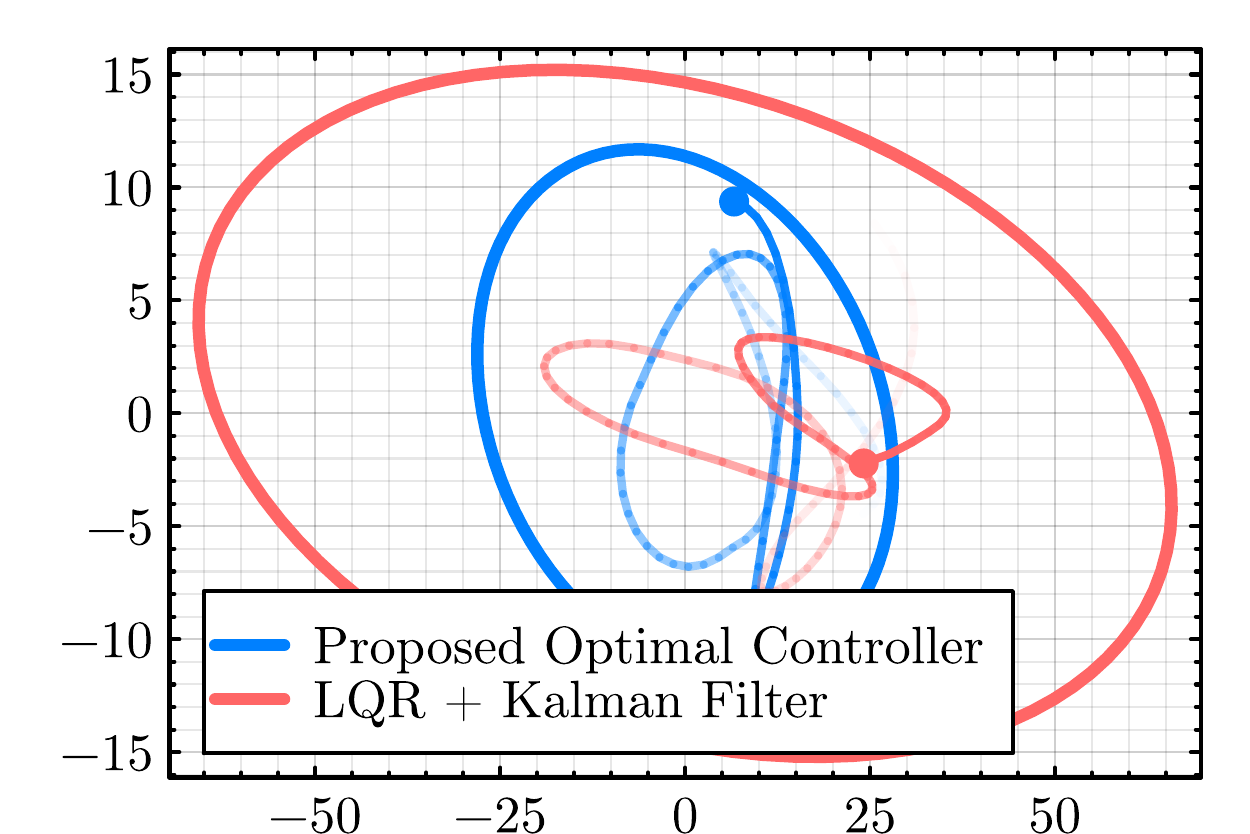}
\makebox[0pt][c]{\hspace{-1.96\textwidth}
    \rotatebox{90}{\begin{minipage}{0.65\textwidth}
        \[ \ z_2 \]
    \end{minipage}}
}\\
\makebox[0pt][c]{
    \begin{minipage}{\textwidth}\vspace*{-0.7cm}
        \[ \ \quad z_1 \]\vspace*{-1cm}
    \end{minipage}
}
\caption{Regulated output $(z_1(t),z_2(t))$ at $t=\frac{T}{3}$}
\end{subfigure}
\hfill
\begin{subfigure}{0.32\textwidth}
\centering
\includegraphics[width=0.98\textwidth]{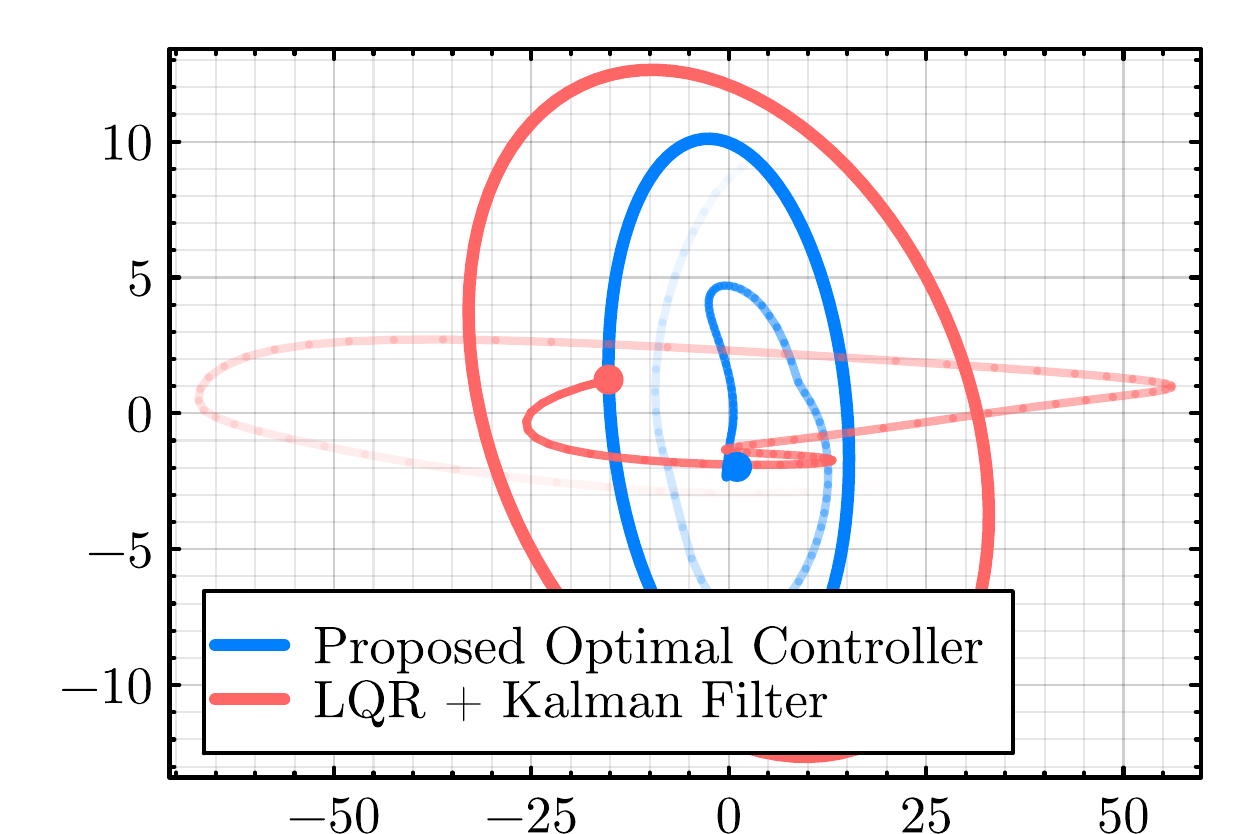}
\makebox[0pt][c]{\hspace{-1.96\textwidth}
    \rotatebox{90}{\begin{minipage}{0.65\textwidth}
        \[ \  z_2 \]
    \end{minipage}}
}\\
\makebox[0pt][c]{
    \begin{minipage}{\textwidth}\vspace*{-0.7cm}
        \[ \ \quad z_1 \]\vspace*{-1cm}
    \end{minipage}
}
\caption{Regulated output $(z_1(t),z_2(t))$ at $t=\frac{2T}{3}$}
\end{subfigure}
\hfill
\begin{subfigure}{0.32\textwidth} 
\centering
\includegraphics[width=0.98\textwidth]{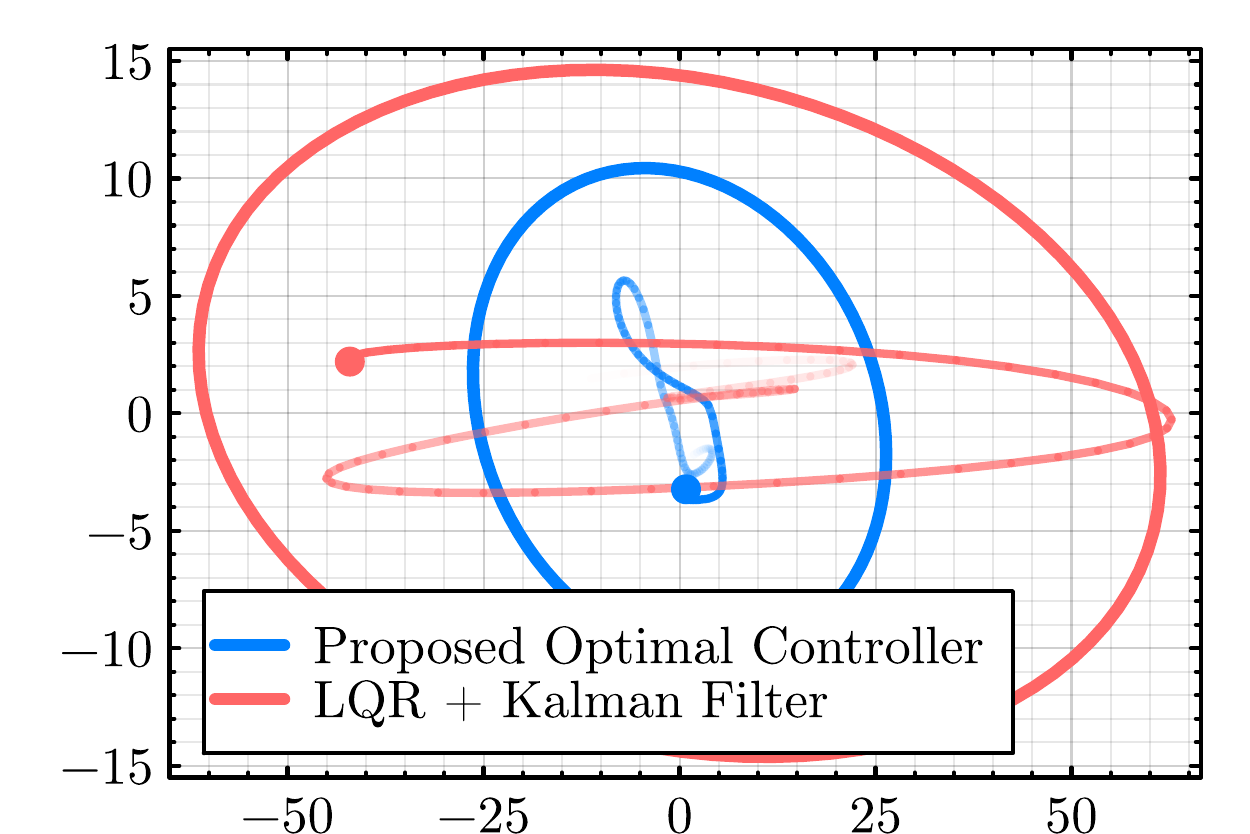}
\makebox[0pt][c]{\hspace{-1.96\textwidth}
    \rotatebox{90}{\begin{minipage}{0.65\textwidth}
        \[ \  z_2 \]
    \end{minipage}}
}\\
\makebox[0pt][c]{
    \begin{minipage}{\textwidth}\vspace*{-0.7cm}
        \[ \ \quad z_1 \]\vspace*{-1cm}
    \end{minipage}
}
\caption{Regulated output $(z_1(t),z_2(t))$ at $t=T$}
\end{subfigure}
\caption{Periodic inescapable ellipsoids in the $z$-plane for the example in Section \ref{Section:Example} with historic trajectories. The smaller ellipsoid is given by the optimal controller \eqref{Output-Feedback} with $\alpha^*(t)$, while the larger one is for the LQR-Kalman filter combination. Note that neither trajectory has left the ellipsoid; the tail outside reflects historical data from when the ellipsoid was larger.}
\label{fig:ellipsoids}
\end{figure*}

Figure \ref{fig:ellipsoids} depicts the periodic inescapable ellipsoids for both the proposed optimal controller and the LQR-Kalman filter combination. The latter results in a larger inescapable ellipsoid, confirming its suboptimality for this criterion, whereas the proposed controller achieves the optimal result.

\section{Conclusion}

This letter presented an approach for optimal control of periodic LTV systems under unknown-but-bounded disturbances. We proposed a systematic method for designing optimal state-feedback, observer, and output-feedback controllers using parameter-dependent Riccati and Lyapunov equations. 
\iffalse
A key contribution was an iterative procedure for optimizing the time-varying parameter $\alpha(t)$, reducing an infinite-dimensional problem to a tractable iterative process.
\fi
Future work may explore robust formulations of the presented results via LMI-based relaxations of the closed-loop design equations.

\appendix

%\section{Proofs of Theorems}
\hypertarget{appendix}{}
\label{appendix}

\begin{prooftheorem}[\ref{Analysis_equation}]
Let $V(x(t),t) \vcentcolon = x(t)^\top \bar{P}(t)x(t)$, where $\bar{P}(t) = P(t)^{-1}$. The ellipsoid $\mathcal{E}_P(t)$ satisfies \eqref{Definition_Inescapability} if and only if, for every trajectory $x(\cdot)$ and every input $w(\cdot)$, the following implication holds at every moment $t$:
\begin{equation} \label{Proof_condition}
    \left\{
    \begin{aligned}
        & V(x(t),t) = 1,\\
        & \lVert w(t)\rVert \le 1,
    \end{aligned}
    \right.
    \quad\Rightarrow\quad
    \dot{V}(x(t),t) \le 0.
\end{equation}
By direct differentiation, expand $\dot{V}(x(t),t)$ as
\begin{align*}
   \dot{V}(x(t),t) &= x(t)^\top \bigl[\dot{\bar{P}}(t) + A(t)^\top\bar{P}(t) + \bar{P}(t)A(t)\bigr]x(t) \\
   &\quad + 2x(t)^\top \bar{P}(t)B(t)w(t).
\end{align*}
Next, note that
\begin{align*}
   2x(t)^\top \bar{P}(t)B(t)w(t) 
   & \overset{(\ast)}{\le} 2 \bigl \lVert B(t)^\top \bar{P}(t)x(t) \bigr \rVert  \\
   & \overset{(\ast\ast)}{\le} \alpha(t) + \frac{ \bigl \lVert B(t)^\top \bar{P}(t)x(t){ \bigr \rVert}^2}{\alpha(t)} \\
   & = \alpha(t) x(t)^\top \bar P(t) x(t)  \\
   & + \frac{1}{\alpha(t)}x(t)^\top \bar{P}(t)B(t)B(t)^\top \bar{P}(t)x(t).
\end{align*}
Both $(\ast)$ and $(\ast\ast)$ are \emph{tight} and hence do not invalidate the necessity and sufficiency of the argument. Indeed, $(\ast)$ becomes an equality under the worst-case input 
\[
    w(t) = \frac{B(t)^\top \bar{P}(t)x(t)}{{ \bigl \lVert}B(t)^\top \bar{P}(t)x(t){ \bigr \rVert}},
\]
and $(\ast\ast)$ is the AM-GM inequality, and there always exists $\alpha(t)$ that turns it into equality. 
% The final equality is justified by the fact that the inequality $\dot{V}(x(t),t) \le 0$ is considered on the boundary $V(x(t),t) =x(t)^\top \bar P(t) x(t) = 1$. 
Therefore, satisfying \eqref{Proof_condition} is equivalent to ensuring that, for some $\alpha(t) >0$,
\begin{align*}
   &\dot{\bar{P}}(t) + A(t)^\top\bar{P}(t) + \bar{P}(t)A(t) \\
   &\quad + \alpha(t)\bar{P}(t) + \frac{1}{\alpha(t)}\bar{P}(t)B(t)B(t)^\top \bar{P}(t) 
   \preceq 0
\end{align*}
Since $\dot{\bar{P}}(t) = -\bar{P}(t)\dot{P}(t)\bar{P}(t)$, multiplying both sides by $P(t)$ yields the desired Lyapunov-type inequality. The minimal $P(t)$ with respect to $\preceq$ makes it an equality (cf. \cite{Abedor1994}).
\endproof
\end{prooftheorem}

\begin{lemma} \label{GR_lemma}
    For $T$-periodic solutions of 
    \begin{align*}
    \dot P(t) &= A(t)P(t) + P(t)A(t)^\top + \alpha(t)P(t) + R(t), \\
    -\dot Q(t) &= A(t)^\top Q(t) + Q(t)A(t)+ \alpha(t)Q(t) + G(t), 
    \end{align*}
    with piecewise-continuous $R(t), G(t) \succeq0$, it holds that
    \[
    \int_0^T \mathrm{trace}(Q(t)R(t))dt = \int_0^T \mathrm{trace}(P(t)G(t))dt.
    \vspace{5px}
    \]
\end{lemma}
\begin{prooflemma}[\ref{GR_lemma}]
    The following identity:
    \[
    \int_0^T \mathrm{trace}(Q(t) \dot P(t))dt
    =
    -\int_0^T \mathrm{trace}(\dot Q(t)  P(t))dt,
    \]
    is obtained using integration by parts. Substituting $\dot P$, $\dot Q$ and rearranging the terms gives the result, where both integrals converge due to the piecewise-continuity of the integrands. %.
    \endproof
\end{prooflemma}

\begin{prooftheorem}[\ref{equivalence_of_norms}]
    Setting $R(t) = \frac{1}{\alpha(t)}B(t)B(t)^\top$, and $G(t) = C(t)^\top  C(t)$ proves the theorem by Lemma \ref{GR_lemma}.
    \endproof
\end{prooftheorem}

\begin{lemma} \label{convexity_lemma}
The size of an ellipsoid $\mathcal{E}_\alpha$ of the system \eqref{System} is strictly convex under variations of $\alpha$.
\end{lemma}
\begin{prooflemma}[\ref{convexity_lemma}]
    The analytical periodic solution of the Lyapunov equation \eqref{P_equation} is:
    \begin{align*}
    P(t) &=
    \int_{-\infty}^{t} \underset{(f \circ g)(\alpha(\tau))}{
    \,%\Bigg[ 
    \underbrace{ \frac{e^{\int_{\tau}^{t} \alpha(s) ds}}{\alpha(\tau)} } 
    \,%\Bigg]
    } \Phi(t,\tau)B(\tau)B(\tau)^{\!\top}\Phi(t,\tau)^{\!\top} d\tau,
    %&=
    %\int_{-\infty}^{t}\mathcal{K}[\alpha](t,\tau) %\Phi(t,\tau)B(\tau)B(\tau)^{\!\top}\Phi(t,\tau)^{\!\top} d\tau,
    %\Phi(t,\tau)B(\tau)B^\top(\tau)\Phi^\top(t,\tau) d\tau,
    \end{align*}
    where $\Phi(t,\tau)$ is the state-transition matrix of the system \eqref{System}.

    Considering that the mapping $f : (x,y) \mapsto \frac{e^y}{x}$ is strictly convex for $x>0$ and the operator $g : \alpha(\tau) \mapsto (\alpha(\tau),\int_\tau^t \alpha(s) ds)$ is linear, their composition $f \circ g$ is strictly convex as well. The mapping $(f \circ g)(\alpha(\cdot)) \mapsto P$ is strictly convex due to controllability, and the mapping $P \mapsto \operatorname{size}(\mathcal{E}_P)$ is strictly convex due to observability.
    \endproof
\end{prooflemma}

\begin{prooftheorem}[\ref{Optimal_alpha_theorem}]
    Following standard optimal control methodology, consider the Hamiltonian associated with the functional $J : \alpha \mapsto \operatorname{size}(\mathcal{E}_\alpha)$:
    \begin{equation*}
        \mathcal{H}(t) \doteq \mathrm{trace} \left( C(t) P(t)C (t)^\top \right) + \mathrm{trace} \left( Q(t) \dot P(t) \right),
    \end{equation*}
    where $Q(t)$ is the Lagrange multiplier, which can be found from the adjoint differential equation $\frac{\partial \mathcal{H}}{\partial P} = -\dot Q(t)$:
    \begin{align*}
         -\dot Q(t) = Q(t)A(t) + A(t)^\top Q(t)  + \alpha(t)Q(t) +C(t)^\top C(t).
    \end{align*}
    The functional derivative of $J$ in the detection of $\alpha$ is
    \[
        \frac{\delta J}{\delta\alpha} = \frac{\partial \mathcal{H}}{\partial \alpha} = \mathrm{trace}\left(Q(t) \left(P(t) - \frac{B(t)B(t)^\top}{\alpha^2(t)} \right) \right).
    \]
    By Lemma \ref{convexity_lemma}, functional $J$ is convex under variations of $\alpha$. Thus, setting the derivative to $0$ yields the unique solution \eqref{Optimal_alpha} proving part (I) of the theorem.

    The terms of functional derivative of $J$ with respect to $\alpha_i$ can be rearranged as follows:
    \[
        \frac{\delta J}{\delta\alpha_i} = \mathrm{trace}(Q_i(t)P_i(t)) \left(1 - \frac{\mathrm{trace}\left(B(t)^\top Q_i(t) B(t) \right)}{\alpha_i^2(t)\mathrm{trace}(Q_i(t)P_i(t))} \right).
    \]
    Substituting the iteration policy \eqref{Fixed_point}, leads to:
    \[
        \frac{\delta J}{\delta\alpha_i} = \frac{\mathrm{trace}(Q_i(t)P_i(t))}{\alpha_i^2(t)} \left(\alpha_i^2(t) - \alpha_{i+1}^2(t) \right).
    \]
    Multiplying both sides by $\Delta \alpha_i = \alpha_{i+1}(t) - \alpha_i(t)$ results in:
    \[
        \frac{\delta J}{\delta\alpha_i}\Delta \alpha_i = -\frac{\mathrm{trace}(Q_i(t)P_i(t))}{\alpha_i^2(t)} \left(\alpha_i(t) + \alpha_{i+1}(t) \right)\Delta \alpha_i^2.
    \]
    Notice that $\frac{\delta J}{\delta\alpha_i}\Delta \alpha_i$ is non-positive for all $t$. Consider the series expansion of the functional $J$:
    \[
    J[\alpha_i+\Delta\alpha_i] = J[\alpha_i] + \int_0^T \frac{\delta J}{\delta\alpha_i} \Delta \alpha_i \,dt + \mu ,
    \]
    where $\mu$ represents higher-order terms that vanish near the optimum. The inequality $J[\alpha_i+\Delta\alpha_i] \leq J[\alpha_i]$ follows, which completes the proof of part (II).
    \endproof
\end{prooftheorem}

\begin{prooftheorem}[\ref{Theorem_KQ}]
    Consider the Hamiltonian function associated with the functional $J : K \mapsto \operatorname{size}(\mathcal{E}_{\alpha,K})$:
    \begin{equation*}
        \mathcal{H}(t) \doteq \mathrm{trace} \left( C_K(t) P(t)C_K (t)^\top \right) + \mathrm{trace} \left( Q(t) \dot P(t) \right).
    \end{equation*}
We derive the optimal controller from $\frac{\partial \mathcal{H}}{\partial K} = 0$:
\begin{align*}
    % &\Rightarrow \frac{\partial \,  \mathrm{trace} \left( K(t) P(t) K^\top(t)\right)}{\partial K} \\
    % &+2 \,\frac{\partial \,  \mathrm{trace} \left( Q(t) B(t) K(t) P(t)  \right)}{\partial K} \, = 0 \\
    P(t)K^\top(t) + P(t) Q(t) B(t) = 0 \; \Leftrightarrow \;
    K(t) = -B^\top(t)Q(t).
\end{align*}
Applying the optimal gain to the adjoint differential equation of system \eqref{State-Feedback_System}-\eqref{State-Feedback_Controller} results in the adjoint Riccati differential equation \eqref{State-Feedback_Q}, which has a unique positive-definite solution as long as the assumptions hold \cite{shayman1984inertia}.%, thereby proving the theorem.
\endproof
\end{prooftheorem}

\begin{prooftheorem}[\ref{Theorem_LP}]
    Consider the Hamiltonian function associated with the functional $J : L \mapsto \operatorname{size}(\mathcal{E}_{\alpha,L})$:
    \begin{equation*}
        \mathcal{H}(t) \doteq \mathrm{trace} \left( C_w(t) P(t)C_w (t)^\top \right) + \mathrm{trace} \left( Q(t) \dot P(t) \right).
    \end{equation*}
We derive the optimal observer from $\frac{\partial \mathcal{H}}{\partial L} = 0$:
\begin{align*}
    C(t) P(t)+ \frac{1}{\alpha(t)} L(t)^\top = 0  \; \Leftrightarrow \;
    L(t) = -\alpha(t) P(t) C^\top(t).
\end{align*}
The rest of the proof is similar to Theorem \ref{Theorem_KQ}.
\endproof
\end{prooftheorem}

\begin{prooftheorem}[\ref{Theorem_KQLP}]
Following the system-splitting approach in \cite[Theorem~8]{peregudin2023} and leveraging the orthogonality of certain systems, one can show that
\begin{equation*} \operatorname{size} (\mathcal{E}_{\alpha, K,L}) = \operatorname{size} (\bar{\mathcal{E}}_{\alpha, K}) + \operatorname{size} (\tilde{\mathcal{E}}_{\alpha, L}), \end{equation*}
where the inescapable ellipsoids $\bar{\mathcal{E}}_{\alpha, K}$ and $\tilde{\mathcal{E}}_{\alpha, L}$ correspond to auxiliary systems, with their sizes determined independently by $K$ and $L$. Consequently, $\operatorname{size} (\mathcal{E}_{\alpha, K,L})$ can be minimized separately over $K$ and $L$. The claim then follows from Theorems \ref{Theorem_KQ} and \ref{Theorem_LP}.
\endproof
\end{prooftheorem}

\bibliographystyle{ieeetr}
\bibliography{file}

\end{document}